\documentclass[11pt]{amsart}
\usepackage{amssymb,amsmath,amsthm}
\usepackage[hmargin=3cm,vmargin=3.5cm]{geometry}
\usepackage{graphicx}
\usepackage{amscd}
\usepackage[all, knot]{xy}
\usepackage{epsfig}
\usepackage{psfrag}
\usepackage{lscape}
\newtheorem{theorem}{Theorem}[section]

\newtheorem{definition}[theorem]{Definition}

\title{Crossing Change Alternating Knots}

\author{SLAVIK JABLAN}

\begin{document}
\maketitle

\begin{abstract}
In this paper we define Crossing Change Alternating Knots (CCA
knots) and their generalization: $k$-CCA knots.
\end{abstract}

\begin{definition}
Let be given a diagram $D$ of a knot (or link). In $D$ we make a
crossing change in every crossing separately, and the rest of the
crossings remain unchanged. From the diagram $D$ with $n$ crossings
we obtain $n$ new diagrams, each with a single crossing changed, and
the corresponding $n$ knots (or links). A diagram $D$ is called {\it
Crossing Change Alternating} (shortly, CCA) if all the knots (links)
obtained by the crossing changes are alternating. A knot (or link)
$K$ is CCA if it has at least one CCA diagram.
\end{definition}

It is clear that a CCA knot (or link) could be alternating, or
non-alternating.

If an alternating knot has a minimal CCA diagram, all its minimal
diagrams are CCA (according to Tait Flyping Theorem). A large class
of CCA knots and links are rational knots and links.

In the case of alternating knots (links) it is sufficient to find
one minimal diagram which is CCA, and all its minimal diagrams will
be CCA. However, this is not true for non-alternating knots (links):
a non-alternating knot (link) can have two different minimal
diagrams, where one is CCA, and the other is not. For example, the
minimal diagram $(2\,1,2)\,(3,-2)$ (Fig. 1a) of the knot $10_{150}$
is not CCA, but its another minimal diagram $8^*-2:.2 0:.-1.-1$ is
CCA.

\begin{figure}[th]
\centerline{\psfig{file=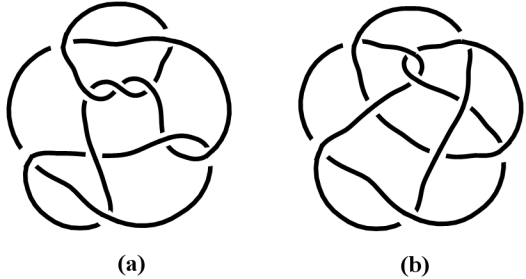,width=3.0in}} \vspace*{8pt}
\caption{(a) The minimal not CCA diagram $(2\,1,2)\,(3,-2)$ of the
knot $10_{150}$; (b) the minimal CCA diagram $8^*-2:.2\,0:.-1.-1$ of
the same knot.\label{fig1}}
\end{figure}

Moreover, CCA-property is not necessarily realized on minimal
diagrams. For example, all minimal diagrams of the knot $10_{151}$
(given in Conway notation as $(2\,1,2)\,(2\,1,-2)$) are not CCA, but
its non-minimal diagram $6^*2\,-1\,-1.2:2\,0$ is CCA, so the knot
$10_{151}$ is a CCA knot without minimal CCA diagrams.

\begin{figure}[th]
\centerline{\psfig{file=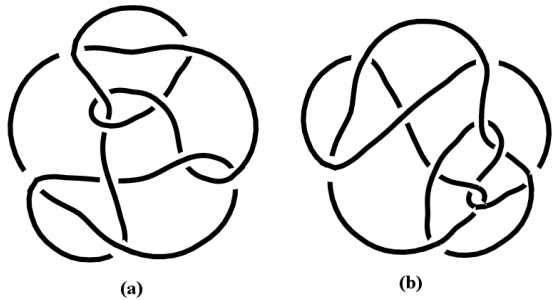,width=3.0in}} \vspace*{8pt}
\caption{(a) The minimal not CCA diagram  $(2\,1,2)\,(2\,1,-2)$ of
the knot $10_{151}$. All its minimal diagrams are not CCA; (b)
non-minimal CCA diagram $6^*2\,-1\,-1.2:2\,0$ of the same
knot.\label{fig2}}
\end{figure}

Thanks to the last example, proving that some knot is CCA is very
difficult, because we need to check all diagrams, and not just the
minimal ones. As the obstruction for a knot to be CCA we can use the
alternation number. The {\it alternation number} of a link $L$,
denoted by $alt(L)$, is the minimal number of crossing changes
needed to deform $L$ into an alternating link, where the minimum is
taken over all diagrams of $L$ \cite{1}. According to T. Abe
\cite{2}, for every knot $K$ alternation number satisfies the
inequality $|{{s(K)-(-\sigma (K))}\over 2}|\le alt (K)$, where
$s(K)$ is the Rasmussen signature of $K$, and $-\sigma (K)$ the
negative signature of $K$. It is clear that any knot $K$ with
$alt(K)>1$ cannot be CCA. T. Abe \cite{2} proved that for every
torus knot $T_{p,q}$ ($2\le p<q$)

\begin{enumerate}
\item $alt (T_{p,q})=0 \Leftrightarrow p=2$;
\item $alt (T_{p,q}=1) \Leftrightarrow (p,q)=(3,4)\,\,or\,\,(3,5)$;
\item $alt (T_{p,q})\ge 2 \Leftrightarrow otherwise$.
\end{enumerate}

Hence, we know that there is an infinite number of knots that are
not CCA. However, the obstruction $|{{s(K)-(-\sigma (K))}\over
2}|\le alt (K)$ is not strong enough for many knots for which we
suspect that are not CCA. E.g., for every alternating knot the
Rasmussen signature and signature coincide, and there are many
alternating knots which are candidates for knots that are not CCA.
Such candidates are all alternating knots with a minimal diagram
which is not CCA.

Definition 1 can be generalized in order to define $k$-CCA knots:

\begin{definition}
Let be given a diagram $D$ of a knot (or link). In $D$ we make $k$
crossing changes in each subset of crossings consisting from $k$
crossings ($1\le k\le \lfloor {n \over 2} \rfloor +1$), and the rest
of the crossings remain unchanged. From the diagram $D$ with $n$
crossings we obtain $n \choose k$ new diagrams, each with $k$
crossings changed, and the corresponding $n \choose k$ knots (or
links). A diagram $D$ is called {\it $k$-Crossing Change
Alternating} (shortly, $k$-CCA) if all the knots (links) obtained by
the crossing changes are alternating. A knot (or link) $K$ is
$k$-CCA if it has at least one $k$-CCA diagram.
\end{definition}

In the same way as before, we expect that there exist knots that are
$k$-CCA, but without a minimal $k$-CCA diagram, so it will be very
difficult to conclude that some knot is $k$-CCA or not. As the
obstruction for a knot to be $k$-CCA we can use the same obstruction
as before. However, based on the computations on minimal diagrams,
there will be many candidates for knots that are not $k$-CCA, for
which will be very difficult to show that they are not $k$-CCA.
Making computations only on the minimal diagrams, we can conclude
that, e.g., the diagram $2\,1,2\,1,2+$ of the knot $9_{28}$ is not
1-CCA, it is 2-CCA, and not 3-, 4-, nor 5-CCA. On the other hand,
the minimal diagram $2\,1,2\,1,-2$ of the knot $8_{20}$ is 1-, 2-,
and 4-CCA, but is not 3-CCA. After checking all minimal diagrams of
this knot, we can conclude that none of them is 3-CCA, but we are
not able to say that the knot $2\,1,2\,1,-2$ is 3-CCA or not,
because we need to check all its non-minimal diagrams, and for this
knot the mentioned obstruction based on alternating number is not
helpful.

\section{CCA knots with $n\le 12$ crossings}

All computations in this paper are made in the program {\it LinKnot}
\cite{3}.

The first table contains alternating CCA knots with $n\le 12$
crossings with minimal CCA diagrams.

\tiny

\bigskip
\noindent
\begin{tabular}{|c|c|c|c|c|c|c|c|c|c|} \hline
$3_1$ &  $4_1$ &  $5_1$ &  $5_2$ &  $6_1$ &  $6_2$ &  $6_3$ &  $7_1$
&  $7_2$ &  $7_3$ \\ \hline $7_4$ &  $7_5$ &  $7_6$ &  $7_7$ &
$8_1$ &  $8_2$ &  $8_3$ &  $8_4$ &  $8_5$ &  $8_6$ \\ \hline $8_7$ &
$8_8$ &  $8_9$ &  $8_{10}$ &  $8_{11}$ &  $8_{12}$ &  $8_{13}$ &
$8_{14}$ &  $8_{15}$ &  $8_{16}$ \\ \hline $8_{17}$ &  $8_{18}$ &
$9_1$ &  $9_2$ &  $9_3$ &  $9_4$ &  $9_5$ &  $9_6$ &  $9_7$ &  $9_8$
\\ \hline $9_9$ &  $9_{10}$ &  $9_{11}$ &  $9_{12}$ &  $9_{13}$ &
$9_{14}$ &  $9_{15}$ &  $9_{17}$ &  $9_{18}$ &  $9_{19}$ \\ \hline
$9_{20}$ &  $9_{21}$ &  $9_{22}$ &  $9_{23}$ &  $9_{25}$ &  $9_{26}$
&  $9_{27}$ &  $9_{29}$ &  $9_{30}$ &  $9_{31}$ \\ \hline $9_{34}$ &
$9_{35}$ &  $9_{36}$ &  $9_{37}$ &  $9_{38}$ &  $9_{39}$ &  $9_{41}$
&  $10_1$ &  $10_2$ &  $10_3$ \\ \hline $10_4$ &  $10_5$ &  $10_6$ &
$10_7$ &  $10_8$ &  $10_9$ &  $10_{10}$ &  $10_{11}$ &  $10_{12}$ &
$10_{13}$ \\ \hline $10_{14}$ &  $10_{15}$ &  $10_{16}$ &  $10_{17}$
&  $10_{18}$ &  $10_{19}$ &  $10_{20}$ &  $10_{21}$ &  $10_{22}$ &
$10_{23}$ \\ \hline $10_{24}$ &  $10_{25}$ &  $10_{26}$ &  $10_{27}$
&  $10_{28}$ &  $10_{29}$ &  $10_{30}$ &  $10_{31}$ &  $10_{32}$ &
$10_{33}$ \\ \hline $10_{34}$ &  $10_{35}$ &  $10_{36}$ &  $10_{37}$
&  $10_{38}$ &  $10_{39}$ &  $10_{40}$ &  $10_{41}$ &  $10_{42}$ &
$10_{43}$ \\ \hline $10_{44}$ &  $10_{45}$ &  $10_{46}$ &  $10_{47}$
&  $10_{50}$ &  $10_{51}$ &  $10_{52}$ &  $10_{53}$ &  $10_{54}$ &
$10_{55}$ \\ \hline $10_{58}$ &  $10_{59}$ &  $10_{60}$ &  $10_{61}$
&  $10_{62}$ &  $10_{63}$ &  $10_{67}$ &  $10_{68}$ &  $10_{69}$ &
$10_{76}$ \\ \hline $10_{77}$ &  $10_{78}$ &  $10_{79}$ &  $10_{82}$
&  $10_{84}$ &  $10_{85}$ &  $10_{87}$ &  $10_{90}$ &  $10_{93}$ &
$10_{100}$ \\ \hline $10_{102}$ &  $10_{103}$ &  $10_{104}$ &
$10_{106}$ &  $10_{108}$ &  $10_{109}$ &  $10_{110}$ &  $10_{111}$ &
$10_{112}$ &  $10_{114}$ \\ \hline $10_{118}$ &  $10_{119}$ &
$10_{120}$ &  $10_{123}$ &  $K11a4$ &  $K11a8$ &  $K11a9$ & $K11a10$
&  $K11a11$ &  $K11a12$ \\ \hline $K11a13$ &  $K11a14$ & $K11a15$ &
$K11a21$ &  $K11a33$ &  $K11a35$ & $K11a37$ &  $K11a39$ &  $K11a42$
&  $K11a45$ \\ \hline $K11a46$ & $K11a49$ &  $K11a50$ & $K11a58$ &
$K11a59$ &  $K11a61$ &  $K11a62$ &  $K11a63$ &  $K11a64$ &  $K11a65$
\\ \hline $K11a74$ &  $K11a75$ & $K11a77$ &  $K11a80$ & $K11a81$ &
$K11a82$ &  $K11a84$ &  $K11a85$ &  $K11a86$ &  $K11a89$
\\ \hline $K11a90$ &  $K11a91$ &  $K11a93$ &  $K11a95$ &  $K11a96$ &
$K11a97$ &  $K11a98$ &  $K11a103$ &  $K11a104$ &  $K11a108$ \\
\hline $K11a110$ &  $K11a111$ &  $K11a117$ &  $K11a119$ &  $K11a120$
&  $K11a121$ &  $K11a123$ &  $K11a133$ &  $K11a134$ &  $K11a135$ \\
\hline $K11a140$ &  $K11a141$ &  $K11a142$ &  $K11a144$ &  $K11a145$
&  $K11a148$ &  $K11a154$ &  $K11a159$ &  $K11a161$ &  $K11a166$ \\
\hline $K11a167$ &  $K11a174$ &  $K11a175$ &  $K11a176$ &  $K11a177$
&  $K11a178$ &  $K11a179$ &  $K11a180$ &  $K11a181$ &  $K11a182$ \\
\hline $K11a183$ &  $K11a184$ &  $K11a185$ &  $K11a186$ &  $K11a188$
&  $K11a190$ &  $K11a191$ &  $K11a192$ &  $K11a193$ &  $K11a195$ \\
\hline $K11a198$ &  $K11a199$ &  $K11a200$ &  $K11a202$ &  $K11a203$
&  $K11a204$ &  $K11a205$ &  $K11a206$ &  $K11a207$ &  $K11a208$ \\
\hline $K11a210$ &  $K11a211$ &  $K11a214$ &  $K11a218$ &  $K11a220$
&  $K11a223$ &  $K11a224$ &  $K11a225$ &  $K11a226$ &  $K11a228$ \\
\hline

$K11a229$ &   $K11a230$ &   $K11a234$ &   $K11a235$ &   $K11a236$ &
$K11a238$ &   $K11a242$ &   $K11a243$ &   $K11a245$ &   $K11a246$ \\
\hline  $K11a247$ &   $K11a249$ &   $K11a250$ &   $K11a256$ &
$K11a258$ &   $K11a259$ &   $K11a260$ &   $K11a263$ &   $K11a268$ &
$K11a278$ \\  \hline  $K11a279$ &   $K11a280$ &   $K11a282$ &
$K11a286$ &   $K11a293$ &   $K11a296$ &   $K11a299$ &   $K11a303$ &
$K11a305$ &   $K11a306$ \\  \hline  $K11a307$ &   $K11a308$ &
$K11a309$ &   $K11a310$ &   $K11a311$ &   $K11a313$ &   $K11a320$ &
$K11a321$ &   $K11a323$ &   $K11a324$ \\  \hline  $K11a325$ &
$K11a330$ &   $K11a333$ &   $K11a334$ &   $K11a335$ &   $K11a336$ &
$K11a337$ &   $K11a339$ &   $K11a341$ &   $K11a342$ \\  \hline
$K11a343$ &   $K11a345$ &   $K11a346$ &   $K11a355$ &   $K11a356$ &
$K11a357$ &   $K11a358$ &   $K11a359$ &   $K11a360$ &   $K11a361$ \\
\hline  $K11a362$ &   $K11a363$ &   $K11a364$ &   $K11a365$ &
$K11a366$ &   $K11a367$ &   $K12a1$ &   $K12a2$ &   $K12a9$ &
$K12a18$ \\  \hline  $K12a20$ &   $K12a22$ &   $K12a24$ & $K12a25$ &
$K12a28$ &   $K12a31$ &   $K12a32$ &   $K12a37$ & $K12a38$ &
$K12a39$ \\  \hline  $K12a54$ &   $K12a55$ & $K12a56$ &   $K12a78$ &
$K12a87$ &   $K12a95$ &   $K12a96$ & $K12a97$ &   $K12a103$ &
$K12a104$ \\  \hline  $K12a105$ & $K12a106$ &   $K12a110$ &
$K12a118$ &   $K12a123$ &   $K12a125$ & $K12a128$ &   $K12a146$ &
$K12a147$ &   $K12a148$ \\  \hline $K12a152$ &   $K12a153$ &
$K12a156$ &   $K12a158$ &   $K12a159$ &
$K12a160$ &   $K12a161$ &   $K12a165$ &   $K12a168$ &   $K12a169$ \\
\hline  $K12a172$ &   $K12a174$ &   $K12a176$ &   $K12a178$ &
$K12a181$ &   $K12a183$ &   $K12a195$ &   $K12a196$ &   $K12a197$ &
$K12a204$ \\  \hline  $K12a206$ &   $K12a216$ &   $K12a217$ &
$K12a221$ &   $K12a226$ &   $K12a229$ &   $K12a238$ &   $K12a239$ &
$K12a241$ &   $K12a243$ \\  \hline  $K12a246$ &   $K12a247$ &
$K12a248$ &   $K12a250$ &   $K12a251$ &   $K12a254$ &   $K12a255$ &
$K12a257$ &   $K12a259$ &   $K12a260$ \\  \hline  $K12a270$ &
$K12a272$ &   $K12a291$ &   $K12a297$ &   $K12a300$ &   $K12a302$ &
$K12a303$ &   $K12a304$ &   $K12a306$ &   $K12a307$ \\  \hline
$K12a327$ &   $K12a330$ &   $K12a331$ &   $K12a344$ &   $K12a345$ &
$K12a353$ &   $K12a356$ &   $K12a357$ &   $K12a358$ &   $K12a360$ \\
\hline  $K12a365$ &   $K12a369$ &   $K12a370$ &   $K12a373$ &
$K12a375$ &   $K12a376$ &   $K12a378$ &   $K12a379$ &   $K12a380$ &
$K12a382$ \\  \hline  $K12a384$ &   $K12a385$ &   $K12a397$ &
$K12a398$ &   $K12a399$ &   $K12a401$ &   $K12a404$ &   $K12a406$ &
$K12a414$ &   $K12a421$ \\  \hline  $K12a422$ &   $K12a423$ &
$K12a424$ &   $K12a425$ &   $K12a433$ &   $K12a436$ &   $K12a437$ &
$K12a443$ &   $K12a444$ &   $K12a447$ \\  \hline  $K12a448$ &
$K12a454$ &   $K12a471$ &   $K12a476$ &   $K12a477$ &   $K12a482$ &
$K12a497$ &   $K12a498$ &   $K12a499$ &   $K12a500$ \\  \hline
$K12a501$ &   $K12a502$ &   $K12a503$ &   $K12a504$ &   $K12a506$ &
$K12a507$ &   $K12a508$ &   $K12a510$ &   $K12a511$ &   $K12a512$ \\
\hline  $K12a514$ &   $K12a515$ &   $K12a517$ &   $K12a518$ &
$K12a519$ &   $K12a520$ &   $K12a521$ &   $K12a522$ &   $K12a526$ &
$K12a527$ \\  \hline  $K12a528$ &   $K12a532$ &   $K12a533$ &
$K12a534$ &   $K12a535$ &   $K12a536$ &   $K12a537$ &   $K12a538$ &
$K12a539$ &   $K12a540$ \\  \hline  $K12a541$ &   $K12a542$ &
$K12a545$ &   $K12a549$ &   $K12a550$ &   $K12a551$ &   $K12a552$ &
$K12a556$ &   $K12a557$ &   $K12a561$ \\  \hline  $K12a563$ &
$K12a564$ &   $K12a565$ &   $K12a568$ &   $K12a569$ &   $K12a573$ &
$K12a576$ &   $K12a577$ &   $K12a579$ &   $K12a580$ \\  \hline
$K12a581$ &   $K12a582$ &   $K12a583$ &   $K12a584$ &   $K12a585$ &
$K12a595$ &   $K12a596$ &   $K12a597$ &   $K12a600$ &   $K12a601$ \\
\hline  $K12a605$ &   $K12a617$ &   $K12a619$ &   $K12a628$ &
$K12a632$ &   $K12a635$ &   $K12a636$ &   $K12a640$ &   $K12a641$ &
$K12a643$ \\  \hline  $K12a644$ &   $K12a646$ &   $K12a648$ &
$K12a649$ &   $K12a650$ &   $K12a651$ &   $K12a652$ &   $K12a653$ &
$K12a657$ &   $K12a663$ \\  \hline  $K12a667$ &   $K12a669$ &
$K12a670$ &   $K12a676$ &   $K12a677$ &   $K12a679$ &   $K12a682$ &
$K12a683$ &   $K12a684$ &   $K12a685$ \\  \hline  $K12a686$ &
$K12a690$ &   $K12a691$ &   $K12a693$ &   $K12a702$ &   $K12a711$ &
$K12a713$ &   $K12a714$ &   $K12a715$ &   $K12a716$ \\  \hline
$K12a717$ &   $K12a718$ &   $K12a719$ &   $K12a720$ &   $K12a721$ &
$K12a722$ &   $K12a723$ &   $K12a724$ &   $K12a725$ &   $K12a726$ \\
\hline  $K12a727$ &   $K12a728$ &   $K12a729$ &   $K12a731$ &
$K12a732$ &   $K12a733$ &   $K12a736$ &   $K12a738$ &   $K12a740$ &
$K12a743$ \\  \hline  $K12a744$ &   $K12a745$ &   $K12a748$ &
$K12a749$ &   $K12a753$ &   $K12a758$ &   $K12a759$ &   $K12a760$ &
$K12a761$ &   $K12a762$ \\  \hline  $K12a763$ &   $K12a764$ &
$K12a770$ &   $K12a773$ &   $K12a774$ &   $K12a775$ &   $K12a784$ &
$K12a786$ &   $K12a789$ &   $K12a790$ \\  \hline  $K12a791$ &
$K12a792$ &   $K12a794$ &   $K12a795$ &   $K12a796$ &   $K12a797$ &
$K12a799$ &   $K12a800$ &   $K12a802$ &   $K12a803$ \\  \hline
$K12a805$ &   $K12a806$ &   $K12a807$ &   $K12a808$ &   $K12a811$ &
$K12a815$ &   $K12a816$ &   $K12a817$ &   $K12a818$ &   $K12a820$ \\
\hline  $K12a824$ &   $K12a826$ &   $K12a827$ &   $K12a833$ &
$K12a834$ &   $K12a838$ &   $K12a839$ &   $K12a840$ &   $K12a842$ &
$K12a843$ \\  \hline  $K12a845$ &   $K12a847$ &   $K12a848$ &
$K12a849$ &   $K12a850$ &   $K12a851$ &   $K12a855$ &   $K12a858$ &
$K12a859$ &   $K12a860$ \\  \hline  $K12a880$ &   $K12a881$ &
$K12a882$ &   $K12a883$ &   $K12a889$ &   $K12a905$ &   $K12a911$ &
$K12a912$ &   $K12a913$ &   $K12a916$ \\  \hline  $K12a920$ &
$K12a929$ &   $K12a937$ &   $K12a938$ &   $K12a950$ &   $K12a952$ &
$K12a955$ &   $K12a969$ &   $K12a970$ &   $K12a972$ \\  \hline
$K12a974$ &   $K12a975$ &   $K12a978$ &   $K12a984$ &   $K12a991$ &
$K12a996$ &   $K12a1015$ &   $K12a1017$ &   $K12a1023$ & $K12a1024$
\\  \hline  $K12a1029$ &   $K12a1030$ &   $K12a1031$ & $K12a1033$ &
$K12a1034$ &   $K12a1039$ &   $K12a1040$ & $K12a1047$ &   $K12a1051$
&   $K12a1052$ \\  \hline  $K12a1060$ & $K12a1063$ &   $K12a1068$ &
$K12a1083$ &   $K12a1089$ & $K12a1093$ &   $K12a1095$ &   $K12a1106$
&   $K12a1107$ & $K12a1114$ \\  \hline  $K12a1115$ &   $K12a1116$ &
$K12a1118$ & $K12a1124$ &   $K12a1125$ &   $K12a1126$ &   $K12a1127$
& $K12a1128$ &   $K12a1129$ &   $K12a1130$ \\  \hline  $K12a1131$ &
$K12a1132$ &   $K12a1133$ &   $K12a1134$ &   $K12a1135$ & $K12a1136$
&   $K12a1138$ &   $K12a1139$ &   $K12a1140$ & $K12a1142$ \\  \hline
$K12a1145$ &   $K12a1146$ &   $K12a1147$ & $K12a1148$ & $K12a1149$ &
$K12a1151$ &   $K12a1153$ & $K12a1156$ &   $K12a1157$ & $K12a1158$
\\  \hline  $K12a1159$ & $K12a1161$ &   $K12a1162$ & $K12a1163$ &
$K12a1164$ & $K12a1165$ &   $K12a1166$ &   $K12a1168$ &   $K12a1169$
& $K12a1171$ \\  \hline  $K12a1176$ &   $K12a1177$ & $K12a1178$ &
$K12a1179$ &   $K12a1180$ &   $K12a1181$ &   $K12a1183$ & $K12a1194$
&   $K12a1202$ &   $K12a1205$ \\  \hline  $K12a1210$ & $K12a1211$ &
$K12a1214$ &   $K12a1220$ &   $K12a1222$ & $K12a1225$ &   $K12a1229$
&   $K12a1240$ &   $K12a1242$ & $K12a1243$ \\  \hline $K12a1244$ &
$K12a1247$ &   $K12a1248$ & $K12a1249$ & $K12a1258$ &   $K12a1259$ &
$K12a1260$ & $K12a1262$ &   $K12a1264$ & $K12a1273$ \\  \hline
$K12a1274$ & $K12a1275$ &   $K12a1276$ & $K12a1277$ &   $K12a1278$ &
$K12a1279$ &   $K12a1281$ &   $K12a1282$ &   $K12a1285$ & $K12a1286$
\\  \hline
$K12a1287$ & $K12a1288$ &    &  &    &  &  &  &  &
\\  \hline

\end{tabular}

\normalsize
\bigskip

The second table contains non-alternating CCA knots with the minimal
CCA diagram.

\small

\bigskip
\noindent
\begin{tabular}{|c|c|c|} \hline
$8_{19}$ & $3,3,-2$ &$\{\{8\},\{6,8,-12,2,14,16,-4,10\}\}$\\
\hline

$8_{20}$ & $3,2\,1,-2$ & $\{\{8\},\{4,8,-12,2,14,16,-6,10\}\}$
\\\hline $8_{21}$ & $2\,1,2\,1,-2$ & $\{\{8\},\{4,8,-12,2,14,-6,
16,10\}\}$\\\hline $9_{42}$ & $2\,2,3,-2$ &
$\{\{9\},\{4,8,18,-14,2,16,-6,10, 12\}\}$\\\hline $9_{43}$ &
$2\,1\,1,3,-2$ & $\{\{9\},\{4,8,10,-14,2,16,18,-6, 12\}\}$\\\hline
$9_{44}$&
$2\,2,2\,1,-2$ & $\{\{9\},\{4,8,-12,2,16,-6,18,10,14\}\}$\\
\hline $9_{45}$ & $2\,1\,1,2\,1,-2$ &
$\{\{9\},\{4,8,10,-16,2,14,18,-6,12\}\}$
\\\hline $9_{46}$ & $3,3,-3$ & $\{\{9\},\{8,-12,16,14,18,-4,-2,
6,10\}\}$\\\hline $9_{47}$ & $8^*-2\,0$ &
$\{\{9\},\{6,8,10,16,14,-18,4,2, -12\}\}$\\\hline $9_{48}$ &
$2\,1,2\,1,-3$ & $\{\{9\},\{4,10,-14,-12,16,2,-6, 18,8\}\}$\\\hline
$9_{49}$ & $-2\,0:-2\,
0:-2\,0$ & $\{\{9\},\{6,-10,-14,12,-16,-2,18,-4,-8\}\}$\\
\hline $10_{150}$ & $6^*.-2.2.2.2\,0$ &
$\{\{10\},\{6,10,16,20,14,2,-18,4,8, -12\}\}$\\\hline $K11n8$ &
$6^*2.2\,1\,0:-3\,0$&
$\{\{11\},\{4,8,16,20,2,-18,6,22,-12,-10,14\}\}$\\
\hline $K11n115$ & $6^*2.-3.2:2\,0$ &
$\{\{11\},\{6,12,16,22,-18,-20,2,8,4, -10,14\}\}$ \\ \hline

$K11n123$ & $6^*-3.2.2.2\,0$&
$\{\{11\},\{6,10,16,22,18,2,-20,8,4,-14,-12\}\}$\\
 \hline

$K11n124$ & $6^*2.-2.2.2.2\,0$ & $\{\{11\},\{6,-10,14,20,-2,18,4,22,
12,8,16\}\}$ \\ \hline

$K11n143$ & $6^*-2.2.-2.2\,0.2\,
0$ & $\{\{11\},\{6,12,-16,22,-18,2,20,-4,-8,14,10\}\}$\\
\hline

$K11n157$ & $9^*-3$ & $\{\{11\},\{6,18,16,12,4,2,-20,-22,10,8,
-14\}\}$ \\ \hline

$K12n147$ & $-2\,-1\,-1,2\,1\,1,2\,1\,1$&
$\{\{12\},\{4,14,18,16,-12,-22,2,24,20,6,-10,-8\}\}$\\
\hline

\end{tabular}

\normalsize
\bigskip

The knot $10_{151}=(2\,1,2)\,(2\,1,-2)$ (Fig. 2) is the example of a
CCA knot without minimal CCA diagram. Its non-minimal diagram
$6^*2\,-1\,-1.2:2\,0$ with the DT code $\{\{11\},\{-4,10,$
$16,20,2,-22,18,8,12,6,14\}\}$ is CCA.

For all the other knots with $n\le 12$ crossings we don't know are
they CCA or not. If they are, they can have only non-minimal CCA
diagrams.

For all knots $K$ up to $n=11$ crossings T. Abe \cite{2} proved that
$alt(K)\le 1$. Hence, we need to consider only knots with $n=12$
crossings.

We computed alternation numbers from all minimal diagrams, and
concluded that all knots $K$ with $n=12$ crossings have the
alternation number $alt(K)\le 1$, except the knots $K12n426$,
$K12n706$, $K12n801$, $K12n835$, $K12n838$, $K12n888$ with the
alternation number 2.

\bigskip

\bigskip

%%%%%%%%%%%%%%%%%%%%%%
%%
%%  REFERENCES
%%
%%%%%%%%%%%%%%%%%%%%%%

\bigskip
\bigskip

\footnotesize

\noindent THE MATHEMATICAL INSTITUTE, KNEZ MIHAILOVA  36, P.O.BOX
367, \\ 11001 BELGRADE, SERBIA

\medskip

\noindent {\it E-mail address:} $\mathrm{sjablan@gmail.com}$

\bigskip

\end{document}